\documentclass[11pt]{article}
\usepackage{amsmath}
\usepackage{amsfonts}
\usepackage{amssymb}
\usepackage{mathrsfs}                  
\usepackage{amsthm}
\usepackage{enumitem}
\usepackage{graphicx}
\usepackage{amscd}
\usepackage{xcolor}
\usepackage[breaklinks,colorlinks,
  citecolor=blue,
  urlcolor=blue]{hyperref}
\usepackage{mathtools}
\usepackage{bbm}
\usepackage{bm}
\usepackage{bbm}
\usepackage[left=1.5cm,right=1.5cm,top=2.5cm,bottom=2cm]{geometry}
\usepackage{caption}
\usepackage{subcaption}
\usepackage[normalem]{ulem}
\theoremstyle{plain}
\newtheorem{theorem}{Theorem}[section]
\newtheorem{proposition}[theorem]{Proposition}

\theoremstyle{definition}
\newtheorem{definition}[theorem]{Definition}

\theoremstyle{remark}

\newcommand{\bfi}{\bfseries\itshape}
\newcommand{\vertiii}[1]{{\left\vert\kern-0.25ex\left\vert\kern-0.25ex\left\vert #1 
    \right\vert\kern-0.25ex\right\vert\kern-0.25ex\right\vert}}

\begin{document}

\title{{\bf Data-driven cold starting of good reservoirs
}}

\author{Lyudmila Grigoryeva$^{1, 2}$, Boumediene Hamzi$^{3, 4, 9}$, Felix P. Kemeth$^{4}$,\\ Yannis Kevrekidis$^{4}$, G Manjunath$^{5}$, Juan-Pablo Ortega$^{6,7}$,\\ and Matthys J. Steynberg$^{8}$}

\date{}
\maketitle

\begin{abstract}
{Using short histories of observations from a dynamical system, a workflow for the post-training initialization of reservoir computing systems is described. This strategy is called cold-starting, and it is based on a map called the starting map, which is determined by an appropriately short history of observations that maps to a unique initial condition in the reservoir space. The time series generated by the reservoir system using that initial state can be used to run the system in autonomous mode, to produce accurate forecasts of the time series under consideration immediately. By utilizing this map, the lengthy ``washouts" that are necessary to initialize reservoir systems can be eliminated, enabling the generation of forecasts using any selection of appropriately short histories of the observations.}
\end{abstract}

\bigskip


\makeatletter
\addtocounter{footnote}{1} \footnotetext{
Universit\"at Sankt Gallen. Faculty of Mathematics and Statistics. Bodanstrasse 6, CH-9000 Sankt Gallen, Switzerland.  {\texttt{Lyudmila.Grigoryeva@unisg.ch} }}
\addtocounter{footnote}{1} \footnotetext{%
Honorary Associate Professor, University of Warwick. Department of Statistics. Coventry CV4 7AL, United Kingdom.
{\texttt{Lyudmila.Grigoryeva@warwick.ac.uk} }}
\addtocounter{footnote}{1} \footnotetext{%
Department of Computing and Mathematical Sciences, Caltech,  Pasadena, CA 91125, US.  
{\texttt{Boumediene.Hamzi@gmail.com} }}
\addtocounter{footnote}{1} \footnotetext{%
Department of Applied Mathematics and Statistics, Johns Hopkins University, Baltimore, MD, US.   
{\texttt{FKemeth1@jh.edu};  \texttt{YannisK@jhu.edu}}}
\addtocounter{footnote}{1} \footnotetext{%
University of Pretoria. Department of Mathematics and Applied Mathematics. Pretoria 0028, South Africa.   {\texttt{Manjunath.Gandhi@up.ac.za} }}
\addtocounter{footnote}{1} \footnotetext{%
Nanyang Technological University. Division of Mathematical Sciences. School of Physical and Mathematical Sciences. Singapore.   {\texttt{Juan-Pablo.Ortega@ntu.edu.sg} } }
\addtocounter{footnote}{1} \footnotetext{%
Honorary Extraordinary Professor, University of Pretoria. Department of Mathematics and Applied Mathematics. Pretoria 0028, South Africa.}
\addtocounter{footnote}{1} \footnotetext{%
University of Pretoria. Department of Physics. Pretoria 0028, South Africa.    {\texttt{Thys.Steynberg@tuks.co.za} }}
\addtocounter{footnote}{1} \footnotetext{%
The Alan Turing Institute, London, UK.}

\makeatother

\tableofcontents

\section{Introduction}\label{Introduction}

{\it Reservoir computing} (RC) \cite{jaeger2001, maass1, Jaeger04, maass2} and in particular {\it echo state networks} (ESNs)  \cite{Matthews:thesis, Matthews1993, Jaeger04} have gained much notoriety in recent years due to their excellent performance in the forecasting of dynamical systems \cite{Jaeger04, pathak:chaos, Pathak:PRL, Ott2018, wikner2021using, arcomano2022hybrid} and to the ease of their implementation. RC aims at approximating nonlinear input/output systems using randomly generated state-space systems (called {\it reservoirs}), in which only a readout map is estimated depending on the learning task. It has been theoretically established that this is indeed possible in a variety of deterministic and stochastic contexts \cite{RC6, RC7, RC8, RC20, RC12}. 

In the context of dynamical systems, it has been shown that this technique has close ties with classical embedding strategies like Takens' Theorem \cite{takensembedding} and generalized synchronizations \cite{kocarev1995general, pecora:synch, ott2002chaos, boccaletti:reports:2002, eroglu2017synchronisation}. See \cite{hart:ESNs, Hart2021, RC18, RC21, manjunath2021universal, berry2022learning} for recent developments in that direction. As we explain in detail later on, this connection implies, in the presence of certain hypotheses, the existence of submanifolds of the state space that are preserved by the reservoir dynamics driven by the observations of the dynamical system that we intend to model. Learning that invariant manifold proves to be beneficial in the dimension reduction of the problem and, more importantly, in the possibility of accurately initializing the reservoir just by using an initial condition of the dynamical system or, alternatively, an appropriately short history of some of its observations. 

This idea has been used for the first time in \cite{kemeth2021initializing} in the context of long-short term memory (LSTM) neural networks, and it is what we call {\it cold-starting} of reservoir systems. Reservoir initialization has traditionally been carried with long washout time series that are used in conjunction with the so called {\it fading memory property} to numerically evaluate the right initial reservoir state. More specifically, there is a collection of conditions that one can impose on the reservoir system to guarantee that, when the length of a time series that is fed into a reservoir tends to infinity, the dependence of the output on the value that was used to initialize the reservoir fades away; see for instance the {\it fading memory property} \cite{Boyd1985}, the {\it echo state property} \cite{jaeger2001, manjunath:prsl, manjunath2022embedding}, or the {\it input forgetting property} \cite{RC9}. Any of these properties imply that if the reservoir is fed with an input for a time sufficiently long, the output that will be obtained will approximate arbitrarily well the state value corresponding to the unique solution consistent with an input defined for all negative times, and all this, we emphasize, regardless of the value that has been used to initialize the reservoir. This input whose only goal is finding an approximating initial state is what we call the ``washout"; the length of the washout necessary for proper initialization depends on the dynamic features of the reservoir that fits the data but it may be quite long, which leads in some instances to a sub-optimal data consumption.

This paper shows that under very general hypotheses, a map can be constructed (we call it the {\it starting map}) that associates with each state of the dynamical system or, equivalently (by Takens' Theorem), a short history of its observations, the unique initial condition in the reservoir space that is consistent with all their past history (the dynamical system is assumed to be invertible). The time series produced by the reservoir system out of that initial condition accurately mimics or path-continues those of the dynamical system that we intend to learn. The availability of this map spares the user from long washout reservoir iterations, which may prove computationally costly and difficult to carry out in the presence of small datasets and, more importantly, allows for immediate prediction. 

The paper is structured as follows. In Section \ref{Good reservoirs and generalized synchronizations} we introduce the notion of good reservoir and we present the reservoir computing forecasting framework in connection with the notion of generalized synchronizations. The reservoir cold-starting methodology is presented in Section \ref{The starting map and cold-starting of reservoir systems} that, as we shall see, is based in the existence of what we call a starting map defined using a synchronization manifold that is obtained as the image of the generalized synchronizations introduced in the previous section. A forecasting method using the starting map is carefully spelled out in this section. Various numerical illustrations that show the pertinence of our methodology are contained in Section \ref{Empirical results}.


\section{Good reservoirs and generalized synchronizations}
\label{Good reservoirs and generalized synchronizations}

In this section we introduce the main tool that will be used in the construction of the starting map described in the introduction, namely, the {\it generalized synchronizations} (GS) between (the observations of) a dynamical system and a reservoir system.

\subsection{Reservoirs and generalized synchronizations} We introduce {\bf reservoir systems} as a state-space system (nonlinear in general) made out of two  equations of the form:
\begin{align} 
\label{eq:RCSystemDet_state}
\mathbf{x}_t &=F(\mathbf{x}_{t-1}, {\bf z}_t), \\
{\bf y}_t&=  h (\mathbf{x} _t),\label{eq:RCSystemDet_eadout}
\end{align}
for all $t \in \mathbb{Z}_- $,  where $F \colon \mathbb{R}^N\times \mathbb{R}^d\longrightarrow  D_N$ and $h: \mathbb{R}^N \longrightarrow {\mathbb{R}}^m $ are the {\bfi reservoir} (randomly generated) and the {\bfi readout} (trainable), respectively. The sequences  ${\bf z} \in (\mathbb{R}^d)^{\mathbb{Z}_-}$ and ${\bf y} \in (\mathbb{R}^m)^{\mathbb{Z}_-}$ stand for the \textbf{input} and the {\bfi  output (target)} of the system, respectively, and $\mathbf{x} \in  (\mathbb{R}^N)^{\mathbb{Z}_-}$ are the associated {\bfi  reservoir states} of dimension $N\in \mathbb{N}^+$, also referred to as the number of virtual {neurons} of the system. In this paper, we are interested in the particular setting where the reservoir \eqref{eq:RCSystemDet_state} is driven by the (in general, partial) observations of a given dynamical system. The learning task consists in the path-continuation of the observations of this dynamical system, or, in a more general and complicated case, in the forecasting of the original dynamical system out of its available partial observations. Hence, for the rest of the paper the inputs and outputs in the system \eqref{eq:RCSystemDet_state}-\eqref{eq:RCSystemDet_eadout} will be chosen according to a particular learning task of interest. Moreover, in both considered learning scenarios only the ability of the reservoir to produce high-precision {\it autonomous} multi-step predictions is assessed.

Let $M$ be a compact finite-dimensional differentiable manifold and let $\phi\in {\rm Diff}^1(M) $ be an invertible discrete-time differentiable dynamical system with differentiable inverse that, for any initial condition $m _0 \in M $,  produces the trajectories $\left\{\phi ^t (m _0)\right\}_{t \in \mathbb{Z}}$. Let $\omega \in C ^1(M, \mathbb{R} ^d) $, $d \in \mathbb{N} $,  be a map that encodes $d$-dimensional observations of the dynamical system and define the $(\phi, \omega)$-{\it delay map} $S _{(\phi, \omega)}:M \longrightarrow \ell^{\infty}(\mathbb{R}^d) $ as $S _{(\phi, \omega)}(m):=\left\{\omega(\phi ^t (m))\right\}_{t \in \mathbb{Z}} $.

Let a reservoir in \eqref{eq:RCSystemDet_state} be a continuously differentiable state map $F: \mathbb{R} ^N\times \mathbb{R} ^d \longrightarrow   \mathbb{R} ^N $ and consider the  drive-response system associated to the inputs ${\bf z}_t = S _{(\phi, \omega)}(m)_t$, $t \in \mathbb{Z}$, that is to the $\omega $-observations of $\phi$ and determined by the recursions:
\begin{equation}
\label{drive-response system}
\mathbf{x} _t=F\left(\mathbf{x} _{t-1}, S _{(\phi, \omega)}(m)_t\right), \quad \mbox{$t \in \mathbb{Z},\, m \in M,$}
\end{equation}
\begin{definition}
    We say that a {\it generalized synchronization} (GS) occurs in this configuration when there exists a map $f:M \longrightarrow \mathbb{R}^N $ (which we call a {\it generalized synchronization})  such that
\begin{equation}
\label{generalized synchronization condition}
 \mathbf{x} _t = f (\phi ^t(m))  \quad \mbox{for any $\mathbf{x} _t  $, $t \in \mathbb{Z} $, and $m \in M $ as in \eqref{drive-response system}.}
\end{equation}
\end{definition}

 The existence of a {\it generalized synchronization} $f$ means that the time evolution of the dynamical system in phase space (not just its observations) drives the response in \eqref{drive-response system}.

\subsection{Good reservoirs} 
\label{Good reservoirs}
The next definition specifies, in terms of the generalized synchronizations that we just introduced, when a reservoir is suitable for the modeling of a given dynamical system. We refer to such systems as {\it good reservoirs}. 

\begin{definition}
\label{good reservoir}
We say that $F: \mathbb{R} ^N\times \mathbb{R} ^d \longrightarrow   \mathbb{R} ^N $ is a {\it good reservoir} for the $\omega $-observations of the dynamical system $\phi\in {\rm Diff}^1(M) $ when it induces a generalized synchronization $f:M \longrightarrow \mathbb{R}^N $ that is also an embedding.
\end{definition}

The term {\it embedding} in the definition means that $f $  is an injective immersion, that is, it is a $C ^1  $ map with injective tangent map and, additionally, the manifold topology in $f(M) $ induced by $f$ coincides with the relative topology inherited from the standard topology in $\mathbb{R}^N $. Equivalently, this means that $f(M) $ is an {\it embedded} submanifold of $\mathbb{R}^N $. 

We emphasize that the existence of generalized synchronizations, in general, and of good reservoirs in particular, is not something generic, and it presupposes that various dynamical constraints are satisfied. We briefly enumerate those constraints and some results in the literature that characterize situations in which they are satisfied. First of all, the definition \eqref{generalized synchronization condition} presupposes that for each $m \in M  $ and the corresponding orbit of observations $S _{(\phi, \omega)}(m) $ there exists a sequence $\mathbf{x} := \left\{\mathbf{x}_t\right\}_{t \in \mathbb{Z}} $ such that \eqref{drive-response system} is satisfied. When that existence property holds and, additionally, the solution sequence $\mathbf{x}  $ is unique, we say that $F$ has the $(\phi, \omega) $-{\it Echo State Property} (ESP) (see \cite{jaeger2001, Manjunath:Jaeger, manjunath:prsl} for in-depth descriptions of this property). Moreover, in the presence of the $(\phi, \omega) $-ESP, the state map $F$ determines a unique causal and time-invariant filter $U ^F: S _{(\phi, \omega)}(M) \longrightarrow (\mathbb{R} ^N)^{\mathbb{Z}} $ that associates to each orbit $S _{(\phi, \omega)}(m) $ the unique solution sequence $\mathbf{x} \in (\mathbb{R} ^N)^{\mathbb{Z}} $ of \eqref{drive-response system}. It can be shown \cite[Lemmas II.2 and II.3]{RC18} that if $F: \mathbb{R} ^N\times  \times \mathbb{R}^d \longrightarrow   \mathbb{R} ^N$ is a continuous reservoir map, then the map 
\begin{equation}
\label{expression synchronization}
\begin{array}{cccc}
f_{(\phi, \omega,F)} :&M &\longrightarrow &\mathbb{R}^N\\
 &m &\longmapsto & p _0 \left(U ^F(S _{(\phi, \omega)} (m))\right), 
\end{array}
\end{equation}
is a generalized synchronization, that is, it satisfies the defining relation \eqref{generalized synchronization condition}. In this expression $p _0:(\mathbb{R} ^N)^{\Bbb Z} \rightarrow \mathbb{R} ^N $ is the projection onto the zero entry of the sequence. More generally, the following relation holds
\begin{equation}
\label{defining claim for GS}
U^F(S _{(\phi, \omega)}(m)) _t=f_{(\phi, \omega,F)} \left(\phi ^t (m)\right), 
\end{equation}
for any $t \in \Bbb Z,\, m \in M$. Additionally, the state synchronization map $f_{(\phi, \omega,F)} $ satisfies the identity:
\begin{equation}
\label{SSM recursion for lemma}
f_{(\phi, \omega,F)} (m)=F \left(f_{(\phi, \omega,F)} (\phi ^{-1}(m)), \omega (m)\right), 
\end{equation}
for all $m \in  M $.

Second, the existence and differentiability of generalized synchronizations need to be addressed. GSs were introduced for the first time in \cite{kocarev1995general}, where it was shown that the asymptotic stability of the system response is a sufficient condition for the existence of a GS. Nevertheless, it was quickly noticed in \cite{pyragas:1996, hunt:ott:1997} that the GS whose existence is guaranteed by this theorem might have poor regularity properties, rendering it useless as an attractor representation and reconstruction tool. This fact motivated the characterization in \cite{hunt:ott:1997} of a first differentiability criterion for GSs. This result has been completed in \cite{RC18} where it was shown that  if $F: \mathbb{R} ^N\times \mathbb{R} ^d \longrightarrow   \mathbb{R} ^N $ is of class $C ^2 $, $\omega $ is of class $C ^1$, and 
\begin{equation}
\label{condition for lfx in differentiable statement}
L_{F _x}< \min \left\{1, 1/ \left\|T \phi ^{-1}\right\|_{\infty}\right\},
\end{equation}
then the map given by \eqref{expression synchronization} is a continuously differentiable GS and it is the only one that satisfies the recursion \eqref{SSM recursion for lemma}. The symbol $L_{F _x} $ in \eqref{condition for lfx in differentiable statement} stands for 
$L_{F _x}=\sup_{(\mathbf{x}, {\bf z})\in \mathbb{R} ^N \times \omega(M)} \left\{ \left\|D _x F(\mathbf{x}, {\bf z})\right\|\right\}$   and $\left\|T \phi\right\|_{\infty}:=\sup_{m \in M} \left\{ \left\|T _m \phi\right\|\right\} $, with $T _m\phi: T _mM \longrightarrow T_{\phi(m)}M $ the tangent map of $\phi$ at $m \in M $. This result is a generalization of the main theorem formulated in \cite{hart:ESNs} for the {\it echo state networks} (ESNs) that we shall introduce later on in \eqref{esn equation}. Moreover, due to the result \cite[Theorem 19]{RC9} and the expression \eqref{expression synchronization}, the synchronization $f_{(\phi, \omega,F)} $ is necessarily Lipschitz with a constant $L_{f_{(\phi, \omega,F)}} $ that satisfies 
\begin{equation}
\label{lipschitz for gs}
L_{f_{(\phi, \omega,F)}}\leq L_{F _x}/\left(1-L_{F _x}\right).
\end{equation}

Finally, there remains the embedding property, which is by far the most elusive of them all when it comes to the formulation of sufficient conditions for it to hold, and that are still not available for very popular reservoir choices like ESNs. To the best of our knowledge, only two general statements are available in this context, both of them for linear reservoirs. The first one is Takens' Theorem \cite{takensembedding, huke:2006} since, in our language, this result shows that in the presence of certain non-resonance conditions and for generic scalar observations $\omega \in C ^2(M, \mathbb{R} ) $ of a dynamical system $\phi\in {\rm Diff}^1(M) $, with $M$ compact and $q$-dimensional, a $(2q+1)$-truncated version $S _{(\phi, \omega)} ^{2q+1}$ of the $(\phi, \omega)$-{delay map} given by  
\begin{equation}
\label{delay map embed}
S _{(\phi, \omega)}^{2q+1}(m):=\left(\omega (m), \omega(\phi ^{-1} (m)), \ldots, \omega(\phi ^{-2q} (m))\right)
\end{equation}
 is a continuously differentiable embedding. This map is in turn the GS corresponding to the linear state map $F(\mathbf{x}, z):=A \mathbf{x}+ \mathbf{C} z $, with $A$ the lower shift matrix in dimension $2q+1  $ and $\mathbf{C}= (1,0, \ldots,0) \in \mathbb{R}^{2q+1} $ which, by Takens' Theorem, constitutes a differentiable GS for the scalar observations of $\phi $. This statement has been generalized in \cite{RC21} where it has been shown that roughly speaking, randomly generated linear systems that have the ESP generate GSs that almost surely have the same properties as Takens' delay embeddings.

\subsection{Good reservoirs are indeed good} 
The next proposition shows that good reservoirs and their associated GS embeddings can be used to adequately represent attractor dynamics in an embedded submanifold of the reservoir space.

\begin{proposition}
\label{conjugation with GS}
Let $F: \mathbb{R} ^N\times \mathbb{R} ^d \longrightarrow   \mathbb{R} ^N $ be  a {\it good reservoir} for the $\omega $-observations of the dynamical system $\phi\in {\rm Diff}^1(M) $ with generalized synchronization $f:M \longrightarrow \mathbb{R}^N $. Then:
\begin{description}
\item [(i)] The set $S:=f(M) \subset \mathbb{R}^N  $ is an embedded submanifold of the reservoir space  $\mathbb{R}^N  $.
\item [(ii)] There exists a differentiable observation map $h: S\longrightarrow \mathbb{R}^d  $ that extracts the one-step-ahead prediction of the observations of the dynamical system out of the reservoir states. That is, with the notation introduced in \eqref{drive-response system} and \eqref{eq:RCSystemDet_eadout}:
\begin{equation}
\label{readout forecast}
h(\mathbf{x} _t)= \omega\left(\phi ^{t+1} (m)\right).
\end{equation}
\item [(iii)]  The maps $F$ and $h$ determine a differentiable dynamical system $\Phi \in C^1(S, S) $ given by 
\begin{equation}
\label{conjugate dynamical system}
\Phi( {\bf s}):=F({\bf s}, h({\bf s})),
\end{equation}
which is $C ^1 $-conjugate to $\phi\in {\rm Diff}^1(M) $ by $f$, that is,
\begin{equation}
\label{c1 conjugation}
f \circ \phi = \Phi \circ f.
\end{equation}
\end{description}
\end{proposition}
 
\noindent\textbf{Proof.\ \ (i)} is an elementary consequence of the fact that $f$ is an embedding (see, for instance, \cite{mta} for details). {\bf (ii)}  
Since the GS $f$ is invertible \textcolor{red}{(on $S$)}, we can consider the map $h:= \omega \circ \phi \circ f ^{-1}: f(M) \subset  \mathbb{R}^N \longrightarrow \mathbb{R} ^d $. Now, using the condition \eqref{generalized synchronization condition}, we have that
\begin{equation*}
h(\mathbf{x} _t)= \omega\circ  \phi(f ^{-1}(\mathbf{x} _t))= \omega\circ \phi \circ \phi ^{t} (m)= \omega\left(\phi ^{t+1} (m)\right),
\end{equation*}
as required. Regarding {\bf (iii)}, it is clear that the map $\Phi  $ defined in \eqref{conjugate dynamical system} is $C^1 $. We now show that it maps into $S$. Let $\Phi({\bf s} ) $ with $ {\bf s} \in S= f(M) $ and let $m \in M  $ such that  ${\bf s} = f (m) $. By the definition of the GS $f$ in \eqref{drive-response system}, we can write ${\bf s}= \mathbf{x} _0 $, where $\mathbf{x} _0 \in \mathbb{R}^N  $ is the zero term of the sequence $\mathbf{x} \in (\mathbb{R} ^N)^{\mathbb{Z}}$ obtained as the output of the system determined by $F$ with the sequence $S _{(\phi, \omega)}(m) $ as input. This implies that
\begin{equation*}
\Phi({\bf s})=F({\bf x}_0, h({\bf x}_0))=F({\bf x}_0, \omega(\phi(m)))= \mathbf{x}_1=f(\phi(m)) \in S,
\end{equation*}
as required. Note that in the second equality, we have used \eqref{readout forecast}, and that the last equality is, once again, a consequence of \eqref{drive-response system}. This equality also proves the conjugation \eqref{c1 conjugation}. \quad $\blacksquare$

\section{The starting map and cold-starting of reservoir systems}
\label{The starting map and cold-starting of reservoir systems}
We now show how the tools that we just introduced can be put to work in the solution of forecasting and path continuation problems for a dynamical system given its observations. The setup of these problems is as follows: suppose that a time series $ \left\{ \omega(m), \omega(\phi(m)), \ldots , \omega(\phi^{T-1}(m))\right\}$ of length $T$ of $\omega$-observations of an invertible dynamical system $\phi\in {\rm Diff}^1(M) $ is provided. In the following paragraphs, we spell out the maps that need to be learned in order to solve the following two problems:
\begin{description}
\item[(i)] The {\bfi path-continuation} at horizon $H \in \mathbb{N}$ of the observations. It consists of determining the values $ \left\{ \omega(\phi^{T}(m)), \omega(\phi^{T+1}(m)), \ldots , \omega(\phi^{T+H-1}(m))\right\}$.

\item[(ii)] The {\bfi  forecasting} of the dynamical system at horizon $H \in \mathbb{N}$. It consists of determining the values $ \left\{ \phi^{T}(m), \phi^{T+1}(m), \ldots , \phi^{T+H-1}(m)\right\}$.
\end{description}
If the functional form of the observation $\omega$ is known, one can obviously obtain a solution for the first problem out of the solution for the second one. 

The solutions to these problems are spelled out in the following theorem in which we assume that we have at our disposal a good reservoir system in the sense of Definition \ref{good reservoir} with generalized synchronization $f:M \longrightarrow \mathbb{R}^N $ and that, moreover, the pair $(\phi, \omega) $ satisfies the necessary conditions for the delay map $S _{(\phi, \omega)}^{2q+1} $ in \eqref{delay map embed}  to be a continuously differentiable embedding via Takens' Theorem, with $q\in \mathbb{N}$ the dimension of $M$. The main ingredient of the following theorem is what we call the {\it starting map} defined as 
\begin{equation}
\label{starting map def}
\sigma:= f \circ \left(S _{(\phi, \omega)} ^{2q+1} \right)^{-1}: S _{(\phi, \omega)} ^{2q+1}(M)\subset \mathbb{R}^{2q+1} \longrightarrow \mathbb{R}^N.
\end{equation}
This terminology is justified by the fact that the starting map produces for each short ($2q+1$)-long history of observations, the unique state value that is consistent with their entire semi-infinite past. Note that if the generalized synchronization $f$ is of the type introduced in \eqref{expression synchronization} and the manifold is compact, then the combination of Takens with the inverse function theorem, together with \eqref{lipschitz for gs} imply that the starting map  $\sigma $ is differentiable and globally Lipschitz.

The proof of the following theorem is a straightforward consequence of Proposition~\ref{conjugation with GS}.

\begin{theorem}[Cold-started forecasting methodology]
\label{theorem starting map}
Let $F: \mathbb{R} ^N\times \mathbb{R} ^d \longrightarrow   \mathbb{R} ^N $ be a good reservoir for the $\omega $-observations of the dynamical system $\phi\in {\rm Diff}^1(M) $. Let $f:M \longrightarrow \mathbb{R}^N $ be the corresponding embedding GS and let $h: S\longrightarrow \mathbb{R}^d$ be the predicting readout introduced in \eqref{readout forecast}. Let  $ \left\{ \omega(m), \omega(\phi(m)), \ldots , \omega(\phi^{T-1}(m))\right\}$ be a sample of $\omega$-observations and assume that $T>2q+1 $. Then:
\begin{description}
\item [(i)] The solution of the {\bfi forecasting problem} is given by the following iterations 
\begin{equation}
\label{iterations for forecasting}
\phi^{T+j}(m)= f ^{-1} \left(F \left(f \left(\phi^{T+j-1}(m)\right), h \left( f \left(\phi^{T+j-1}(m)\right) \right) \right)  \right),\quad \mbox{$j=0, \ldots, H-1$,}
\end{equation}
that can be readily initialized at $j=0 $ if the state $\phi^{T-1}(m) $ is known. If only observations are available, then the {\it starting map} $\sigma: \mathbb{R}^{2q+1} \longrightarrow S $ defined as $\sigma:= f \circ \left(S _{(\phi, \omega)} ^{2q+1} \right)^{-1} $  has to be used and applied to a $(2q+1) $-long history of observations preceding the instant $T-1 $, which yields:
\begin{equation}
\label{starting map}
\sigma \left(\omega(\phi^{T-1}(m)), \omega(\phi^{T-2}(m)), \ldots , \omega(\phi^{T-2q-1}(m)) \right)=f \left(\phi^{T-1}(m) \right)
\end{equation}
and can be used to initialize the iterations \eqref{iterations for forecasting} at $j=0 $.
\item [(ii)] The solution of the {\bfi path-continuation problem} is given by the following iterations 
\begin{align}
\label{path continuation states}
\mathbf{x}_{T+j-1}&=F \left(\mathbf{x}_{T+j-2}, \omega\left( \phi^{T+j-1}(m)\right) \right),\\
\omega\left( \phi^{T+j}(m)\right)&=h \left(\mathbf{x}_{T+j-1}\right), \quad j=0, \ldots, H-1,\label{path continuation prediction}
\end{align}
where \eqref{path continuation states} is initialized by setting 
\begin{equation}
\label{starting map path}
\mathbf{x}_{T-2}= f \left(\phi^{T-2}(m)\right)= \sigma \left(\omega(\phi^{T-2}(m)), \omega(\phi^{T-3}(m)), \ldots , \omega(\phi^{T-2q-2}(m)) \right).
\end{equation}
\end{description}
\end{theorem}

\subsection{The forecasting method and implementation} 
\label{The forecasting method and implementation}
The forecasting approach contained in Proposition \ref{conjugation with GS} and in Theorem \ref{theorem starting map} requires a few ingredients. More explicitly, first, one needs to devise a good reservoir $F: \mathbb{R} ^N\times \mathbb{R} ^d \longrightarrow   \mathbb{R} ^N $ for the $\omega$-observations of the dynamical system $\phi \in \text{Diff}^1(M)$, $\text{dim}(M) = q$, under consideration. Second, a predicting readout map $h:\mathbb{R}^N \rightarrow \mathbb{R}^d$ introduced in \eqref{readout forecast} needs to be constructed. Finally,  the corresponding GS $f:M\rightarrow \mathbb{R}^N$ and a starting map $\sigma:\mathbb{R}^{2q+1}\rightarrow \mathbb{R}^N$, for initializing the states of the reservoir $F$ in order to construct autonomous multi-step predictions out of the short, $(2q+1)$-long, histories of the dynamical system's observations, need to be obtained. In the following paragraphs we spell out the details of the choice of the design for our forecasting experiment in the next section.

\medskip

\noindent{\bf The reservoir:} We shall be using a leaking {\it echo state network} (ESN) given by
\begin{equation}
\label{esn equation}
F(\mathbf{x} , {\bf z}):=(1- \alpha)\mathbf{x}+ \alpha \tanh  \left(A\mathbf{x}+ C {\bf z} \right),
\end{equation}
where $\alpha \in (0,1]$ is a prespecified {\it leak rate}, $A$ is a square randomly generated {\it connectivity matrix} of dimension $N\in \mathbb{N}$, and $C $ is an {\it input matrix} of dimension $N \times d $ that connects the $d $-dimensional $\omega  $-observations of the dynamical system $ \phi $ to the reservoir given by $F$. The random parameters are sampled such that the sufficient condition $\|A\|_2<1$ ($\|\cdot\|_2$ denotes the matrix 2-norm) for the $(\phi, \omega) $-Echo State Property to hold (see Subsection~\ref{Good reservoirs}) is  satisfied (see, for example, \cite{RC9}). In practice, since the spectral radius $\rho(A)$ satisfies that $\rho(A)\leq \|A\|_2$ it suffices to take $\rho(A)<1$, which is the most common condition used in the reservoir computing literature.

\medskip

\noindent{\bf The forecasting readout:} ESNs have been shown to be universal input-output approximants with linear readouts \cite{RC7, RC20}. This implies that one can choose the predicting readout $h$ to be a linear map, though any choice of a higher-order polynomial or a neural network function is also possible.  A geometric intuition behind the possibility of achieving universal approximation using exclusively linear readouts in reservoir computing has been provided in \cite{RC13, RC19}. Those references show that some universal reservoir computing families that use linear readouts (the so-called state-affine systems (SAS) \cite{RC6}, in this case) are random projections of Volterra series expansions with semi-infinite inputs. Volterra series are an infinite dimensional object whose universality has been proved in \cite{Boyd1985}, and the Johnson-Lindenstrauss Lemma \cite{JLlemma} can be used to show that universality is preserved under the random projections which yield (universal) SAS. We emphasize that this argument applies exclusively to SAS. An analogous result for ESNs remains an open problem.

Later on in Section~\ref{Empirical results} we shall be presenting results that are obtained assuming that the forecasting readout $h$ introduced in \eqref{readout forecast} has a linear functional form, that is, $h( \mathbf{x})= {W}\mathbf{x}  $, where  ${W} $ is a $ d \times N $ matrix. The $T$-long sample of (partial) $d$-dimensional observations of a given dynamical system is used to drive the reservoir $F$ starting from the initial state $\mathbf{x}_0\in \mathbb{R}^N$ chosen to be either a zero vector or a randomly sampled vector. The corresponding $T$ states are collected during this phase, sometimes called the {\it listening phase} in the literature \cite{Verzelli2020b}. To eliminate the influence of the original initialization, the readout map is estimated after discarding the first $T_w$ observations, which is sometimes called the {\it washout period}. This is the most popular approach in the successful applications of reservoir computing cited in the introduction. 

We point out that when we are working on a path-continuation problem dealing with low-dimensional observations of the dynamical system, we shall most likely be agnostic with respect to the dimension $q$ of the data-generating dynamical system $\phi $. This is a classical and well-studied problem that appears when using embedding techniques in dynamical systems forecasting; some techniques for the estimation of the dimension $q$ can be found in \cite{kantz:Sreiber, martin2019impact} and references therein. 

Subsection~\ref{Robustness in learning the cold-start map} contains results obtained without restricting the readouts of the ESNs to be exclusively linear. Instead, a neural network $h_{NN}$ of relatively simple architecture is used as a readout. This renders the forecasts to be derived as a function $h_{NN}$ of the iterates of the system in the autonomous run defined by  $\mathbf{x}_{n+1} =  F(\mathbf{x}_n, h_{NN}(\mathbf{x}_n))$, where the starting map $\sigma$ is applied only once to obtain the initial condition (for instance using $\mathbf{x}_{n} = \sigma \left(\omega(\phi^{n}(m)), \omega(\phi^{n-1}(m)), \ldots , \omega(\phi^{n-2q}(m)) \right)$).  By opting for a neural network as a readout, we can more reliably investigate the potential detrimental effects of an imprecise estimation of the starting map on the accuracy of ESN predictions, rather than including the linear approximation of the readout also contributing to such effects.


\medskip

\noindent{\bf The synchronization manifold $S $ and the starting map $\sigma$.} If the reservoir devised in the previous points is good (in the sense of Definition \ref{good reservoir}), then it has an associated  GS map $f:M \longrightarrow \mathbb{R}^N$  whose image $S=f(M)$ is an embedded submanifold that is left invariant by the reservoir dynamics. Moreover, as we saw in Theorem \ref{theorem starting map}, there exists a map $\sigma: \mathbb{R}^{2q+1} \longrightarrow S $ that links short histories of the dynamical system observations to the points in the synchronization manifold $S $ that are the images of the unique points in phase space $M$ that these histories represent according to Takens' theorem. In Section~\ref{Empirical results} we adopt two techniques for the learning of the starting map, namely, (i) a diffusion maps-based methodology \cite{coifman2006diffusion} which allows to learn together with the starting map the associated synchronization manifold out of the data, and (ii) static feed-forward neural networks. These objects are known to be dense in the set of continuous functions with respect to the topology of uniform convergence \cite{cybenko}, which, in particular, guarantees the learnability of the starting map $\sigma$ since, as we already pointed out in the discussion after \eqref{starting map def}, this map is differentiable under very general conditions.

\medskip

\noindent{\bf The forecasted and path-continued values.} They are obtained by using the recursions and the initializations spelled out in \eqref{iterations for forecasting} and in \eqref{path continuation states}-\eqref{path continuation prediction}, respectively. We note that, unlike in the path-continuation problem, the solution of the forecasting problem requires the learning of not only the synchronization map $f$ but also of its inverse $f ^{-1}$. We hence restrict our empirical analysis in Section~\ref{Empirical results} to the case of the path-continuation learning problem. 

\medskip

\noindent{\bf Importance of the informed cold-starting.} The most important difference between the methodology that we just proposed and the one used in all the above-cited empirical contributions is in the reservoir initializations proposed in the equations \eqref{starting map} and \eqref{starting map path} for the forecasting and path-continuation problems, respectively. More explicitly, having obtained the readout map using some chosen loss function and solving the associated empirical risk minimization (ERM) problem (for example, for linear readouts and quadratic losses the solutions of the corresponding ERM problems are the least squares solutions), one would traditionally reason as follows: given a history of observations for the path-continuation and the forecasting problems, one needs to initialize the reservoir state to construct the predictions. In the traditional approach, the initialization values $\mathbf{x}_{T-1}=f \left(\phi^{T-1}(m) \right) $  and  $\mathbf{x}_{T-2}=f \left(\phi^{T-2}(m) \right) $ in \eqref{starting map} and \eqref{starting map path}, respectively, are obtained by feeding a  sequence of observations $\left\{\omega(\phi^{T-n}(m)), \ldots, \omega(\phi^{T-1}(m))\right\} $ into the reservoir that is initialized at an arbitrary state $\mathbf{x} _0 \in \mathbb{R}^N $.  Subsequently, the last state is processed with the trained readout map and the output is used to autonomously run the reservoir for the desired number of future steps of the multi-step path-continuation or forecasting exercise. It is well known that for short history sample of observations used as inputs this traditional approach would lead to poor predicting performance of the reservoir since, in this case, the impact of the initialization of the states is very high. More explicitly, consider the iterations
\begin{equation*}
\mathbf{x} _{T-j}^n(\mathbf{x} _0)=F\left(\mathbf{x} _{T-j-1}^n(\mathbf{x} _0), \omega(\phi^{T-j-1})\right), \quad \mbox{$j \in \left\{1, \ldots, n\right\},\  \mathbf{x}^n_{T-n-1}= \mathbf{x} _0 \in \mathbb{R}^N.$}
\end{equation*}
Systems that are traditionally used in RC have the so-called fading memory property \cite{Boyd1985}, and, in particular, the input forgetting property \cite{RC9}, which implies that:
\begin{equation*}
\lim _{n \rightarrow \infty} \mathbf{x} _{T-1}^n(\mathbf{x} _0)=\mathbf{x} _{T-1}=f \left(\phi^{T-1}(m) \right), \quad \mbox{for any $\mathbf{x} _0 \in \mathbb{R}^N$.} 
\end{equation*}

We find that our approach offers significant improvements compared to {\it traditional modus operandi}. More precisely, initializing the reservoir with the image of the learned starting map $\sigma$ and hence ``informing'' the original state of the reservoir about the commencing point of our forecasting exercise leads to less data-intensive predictions since no washout periods are needed. Using short histories of observations of length $2q+2 $ for the path-continuation problem and $2q+1$ for the forecasting problem, we can immediately work out what the next time series value is just by using the iterations \eqref{iterations for forecasting} or \eqref{path continuation states}-\eqref{path continuation prediction}. 
The cold-starting procedure that we propose in Theorem \ref{theorem starting map} based on learning the starting map $\sigma $ circumvents the asymptotic traditional approach that may prove costly both from the computational and the data consumption points of view and does not allow to produce high-quality multi-step predictions based on a data of limited length ($2q+2 $ and $2q+1$ for the path-continuation and the forecasting problem, respectively).

\section{Empirical results}
\label{Empirical results}

\begin{figure}[h!]
	\centering
	\begin{subfigure}{0.4\textwidth}
        \centering  \includegraphics[width=0.75\linewidth]{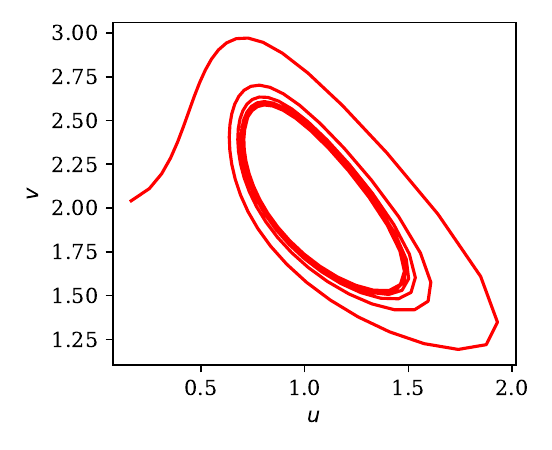}
        \caption{}
    \end{subfigure}
    \begin{subfigure}{0.4\textwidth}
        \centering
\includegraphics[width=0.75\linewidth]{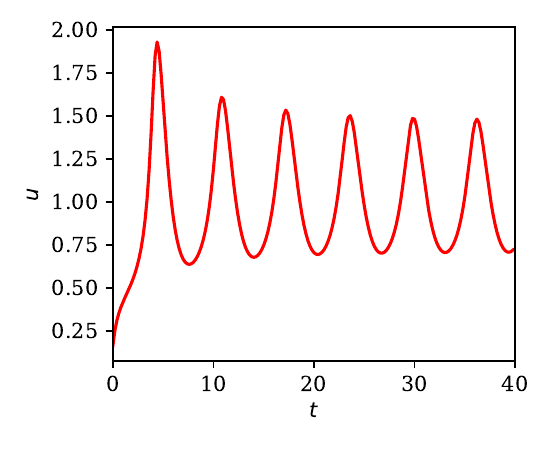}
\caption{}
\end{subfigure}
\caption{\footnotesize  Representative trajectory of the Brusselator system sampled with $\delta t=0.2$. Initial conditions are drawn uniformly such that $u_0\sim \mathcal{U}[0,2]$ and  $v_0\sim \mathcal{U}[0,3]$. (a) Trajectory in phase space of the Brusselator system. (b) $u$ variable evolution of trajectory in (a).}
\label{fig:brusselator_data}
\end{figure}

In this section, we demonstrate the empirical forecasting improvements exhibited by our proposed cold-starting of the reservoir with respect to traditional approaches. We shall use two dynamical systems, namely, the Brusselator, and the Lorenz systems.
The Brusselator is a two-dimensional ($q=2$) system exhibiting oscillatory dynamics~\cite{kondepudi14_brusselator_chapter} given by 
\begin{align*}
    \dot{u} &= a + u^2v - (b+1)u,\\
    \dot{v} &= bu - u^2v,
\end{align*}
and parametrized by $a=1$ and $b=2.1$. For this set of parameters $a$ and $b$,
the only stable attractor of the Brusselator is a stable limit cycle. Figure~\ref{fig:brusselator_data} provides a representative trajectory in phase space and the temporal evolution of $u$ over time.

\begin{figure}[h]
	\centering
	\begin{subfigure}{0.4\textwidth}
        \centering  \includegraphics[width=0.75\linewidth]{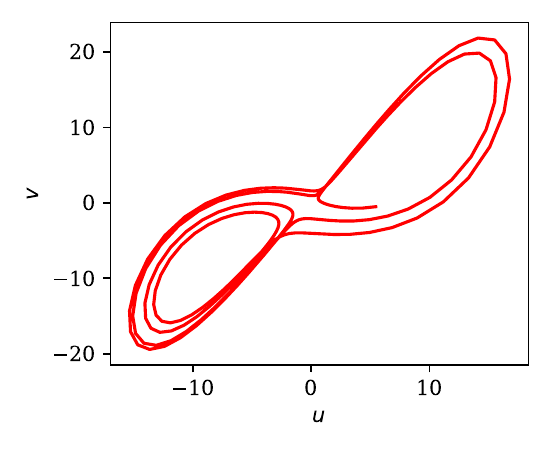}
        \caption{}
    \end{subfigure}
    \begin{subfigure}{0.4\textwidth}
        \centering
\includegraphics[width=0.75\linewidth]{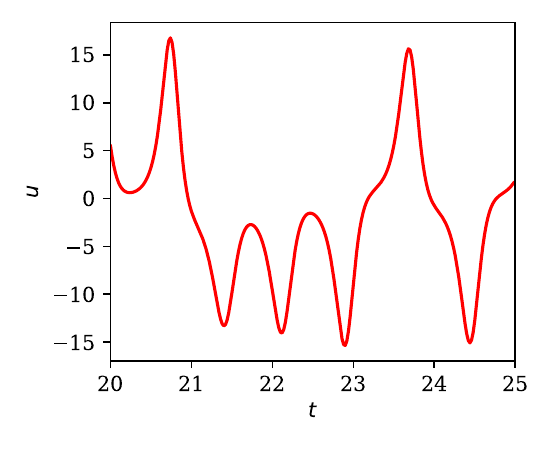}
\caption{}
\end{subfigure}
\caption{\footnotesize  Representative trajectory of the Lorenz system sampled with $\delta t=0.2$. Initial conditions $u_0\sim \mathcal{N}(10,1)$, $v_0\sim \mathcal{N}(1,1)$, and $w_0 \sim \mathcal{N}(0,1)$. (a) Projection of this trajectory onto the $u$-$v$ plane. (b) $u$ variable evolution of trajectory in (a).}
\label{fig:lorenz_data}
\end{figure}

The Lorenz system is a dynamical system presenting a simplified three-dimensional model ($q=3$) for weather prediction \cite{lorenz1963deterministic} and is given by
\begin{equation*}
\begin{aligned}
& \dot{u}=a \left(v-u\right), \\
& \dot{v}=b u-u w-v, \\
& \dot{w}=u v-c w,
\end{aligned}
\end{equation*}
where we use the parameters $a=10$, $b=28$, and $c=8/3$.
For this set of parameters, the dynamics of the Lorenz system exhibits chaotic motion. In Figure~\ref{fig:lorenz_data} a projection of the phase space on the $u$-$v$ plane and the temporal evolution of $u$ are provided.

For these systems, we assume that only the first coordinate observations are available for the learning, and we are interested in their path continuation for $H$ steps into the future based on $2q+1$ past observations. The ESN reservoir systems as in  \eqref{esn equation} are implemented and the forecasting method discussed in Subsection~\ref{The forecasting method and implementation} is followed for the path continuation exercise. To illustrate the forecasting performance, the readout map $h$ is assumed to be linear. For each system, $T$ input  observations of the first coordinate ($d=1$) discretized at $\delta t$ are used to run the estimation procedure. More precisely, after discarding the first  $T_w$-long washout of the states and denoting $T_{\text{tr}}:=T- T_w - 1$, we construct $X:=(\mathbf{x}_{T_w+1}|\mathbf{x}_{T_w+2}|\cdots|\mathbf{x}_{T_w + T_{\text{tr}}} ) \in \mathbb{R}^{N \times T_{\text{tr}}}$ using demeaned states, as well as $U:=({u}_{T_w+2}|u_{T_w+3}|\cdots|u_{T_w + T_{\text{tr}} + 1} ) $ using demeaned one-step ahead true observations of the coordinate $u$. The estimated linear readout map $\mathbf{W}\in \mathbb{R}^{N\times 1}$ is given by the following closed-form solution of the ridge regression:
\begin{equation}
\label{OLS_reservoir}
  \widehat{\mathbf{W}}_{ridge} = \left(X X^\top +\lambda \mathbb{I}_N\right)^{-1}X U^\top,
\end{equation}
with the ridge regularization penalty  $\lambda>0$. The estimated readout is hence  defined by $\widehat{h}_{ridge}(\mathbf{x}) =  \hat{\mathbf{W}}_{ridge} \mathbf{x}$, $\mathbf{x}\in \mathbb{R}^N$. Once the readout map $\widehat{h}_{ridge}$ is available, we compare the performance of the autonomous multi-step  path-continuation of the first coordinate history for each of the systems adopting {\bf{(i)}} the traditional way of initializing the states of the reservoir for the predicting exercise, {\bf{(ii)}} using the starting map proposed in this paper.
For the Brusselator system, $600$ trajectories with initial conditions $u_0\sim \mathcal{U}[0,2]$ and  $v_0\sim \mathcal{U}[0,3]$ are sampled with $\delta t=0.2$ for $30$ dimensionless time units, which results in $T=150$. We discard the first $T_w=1$ washout discretized steps. This results in $600$ pairs of $T_{\text{tr}}=148$-long training paths. For the testing phase, we create trajectories with initial conditions drawn from the same uniform distribution but recorded for $40$ dimensionless time units (200 discrete steps). One testing trajectory is depicted in Figure~\ref{fig:brusselator_data}.
We implement the reservoir with $N=1024$ and $\alpha=0.51$, while taking  $\lambda=0.01$ in \eqref{OLS_reservoir}. The spectral radius of the reservoir matrix is taken as $0.98$. All reservoir hyperparameters are obtained with by performing hyperparameter optimization using the Optuna framework~\cite{optuna_2019}. 

\begin{figure}[h!]
  \centering
  \includegraphics[width=0.7\textwidth]{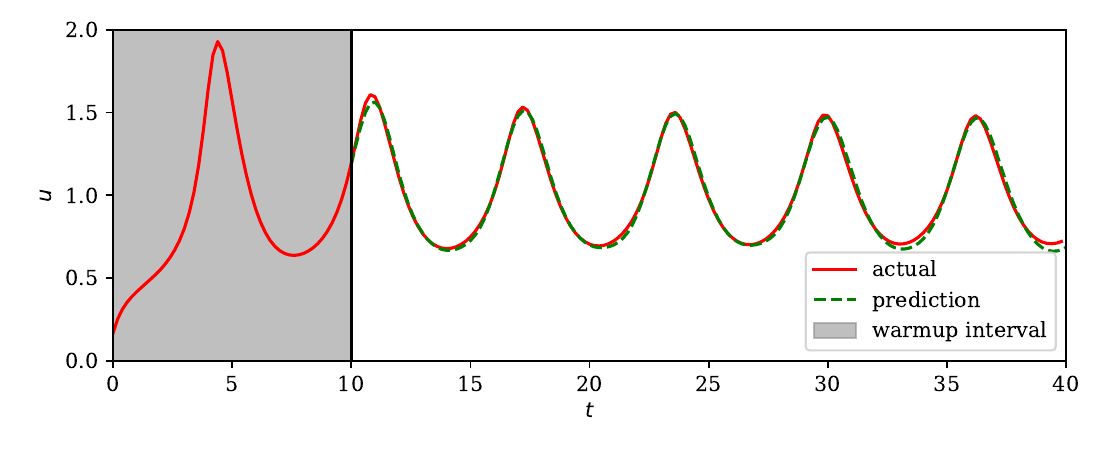}
  \caption{{\footnotesize Autonomous path-continuing of the partial observations of the Brusselator system produced by the ESN with the readout $\widehat{h}_{ridge}$ and with the initial zero state compared to the true trajectory. The shaded area marks the part of the path which is used as a history to drive the trained ESN and the black line shows the moment when the subsequent autonomous path-continuation starts.}}
  \label{fig:prediction}
\end{figure}

To showcase the traditional approach {\bf{(i)}}, we attempt at path-continuing the first coordinate's observations of the Brusselator system for $H=150$  steps into the future (this corresponds to the $30$ steps in the system's time). The standard \textit{modus operandi} consists in initializing the states with any arbitrary starting value, for example, with a randomly sampled or a zero vector (as in our case), 
force the reservoir with a warmup trajectory,
collect the last state corresponding to the last observation of the history of observations that needs to be continued and thereby autonomously iterate the reservoir system with the prior trained readout map to produce $H$ forecasts.
Figure~\ref{fig:prediction} shows an example of this approach for the Brusselator system.
There, the trained reservoir is warmed up by providing an initial warmup $u$ trajectory of length $50$ steps ($10$ dimensionless time units) as input to the model, indicated by the gray-shaded region.
We take the warmup long enough as to approximately washout the influence of the initialization.
The ESN is then used in an autoregressive fashion for $H=150$ steps ($30$ dimensionless time units), producing forecasts for the initial input time series.
The produced forecasts (dashed green curve) are shown together with the actual dynamics of $u$ (solid red curve).
The figure thus illustrates that the trained reservoir model is able to accurately continue the dynamics. In addition, a warmup length of only $50$ time steps ($10$ dimensionless time units) is sufficient to synchronize the internal reservoir states to the input trajectory.

\begin{figure}[h!]
  \centering
  \includegraphics[width=0.7\textwidth]{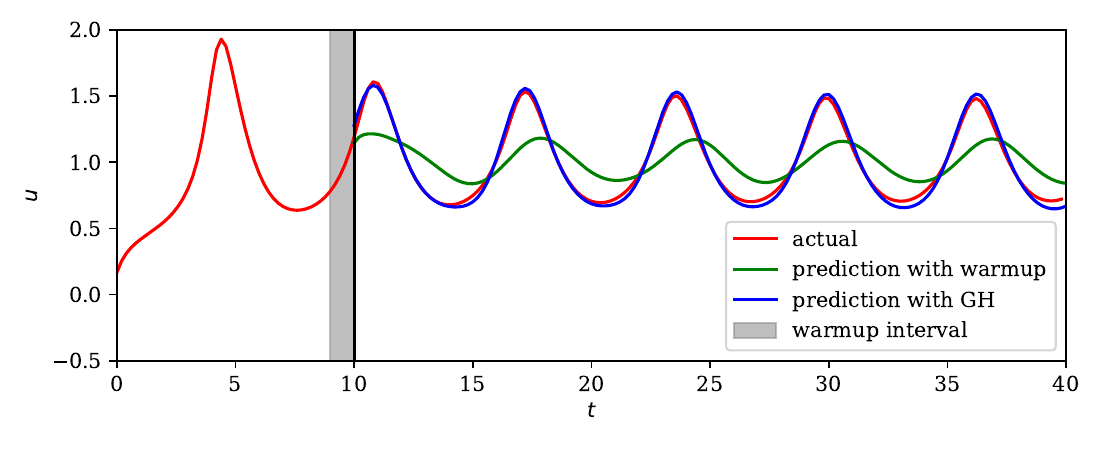}
  \caption{ \footnotesize Representative trajectory of the test data for the Brusselator system (red). The warmup period is of length 5 (gray-shaded region). Green - predictions of the ESN with the states initialized as zero vectors and warmup used.
  Blue - predictions of the ESN with the states initialized with geometric harmonics (GH)  method.}
  \label{fig:predictions_new}
\end{figure}

We emphasize that, in general, our goal is to be able to produce accurate path-continuation using only a minimally short history of these observations, for example, $2q+1=5$, which is possible with our cold-starting technique. For the traditional approach, though, $5$ observations are insufficient to remove the influence of the arbitrary starting initialization. Figure~\ref{fig:predictions_new} demonstrates this scenario. 

We now refer to our cold-starting technique using our proposed starting map $\sigma$. 
Following the same approach as suggested for LSTM networks in \cite{kemeth2021initializing}, we apply diffusion maps to input time series windows (here of length $5$, sampled from the training trajectories) to learn the data manifold as a first step~\cite{Lehmberg2020}.
The two independent diffusion modes $\boldsymbol{v}^{(1)}$ and $\boldsymbol{v}^{(2)}$ that span the data manifold are depicted in Figure~\ref{fig:brusselator_embedding_space} (see Appendix~\ref{sec:dmaps} for the detailed calculation of these modes). Note that each dot corresponds to a time series window of $u$ of length $5$. In addition, notice that the two-dimensional embedding obtained this way is in agreement with the dimensionality of the Brusselator system ($q=2$).
More precisely, for each of the training trajectories, we also produce trajectories of forced internal reservoir states.
We thus obtain for each time series window also corresponding (approximately warmed-up) internal states $\mathbf{x}_i$.
In Figure~\ref{fig:brusselator_embedding_space}, we color each window with one hidden state variable $\mathbf{x}_0$ that corresponds to the last time step of each window.
We now learn a mapping from the (two-dimensional) data manifold to warmed-up internal states of the reservoir using geometric harmonics as it is done for the case of LSTM recurrent neural networks in \cite{kemeth2021initializing}.

\begin{figure}[ht]
  \centering
  \includegraphics{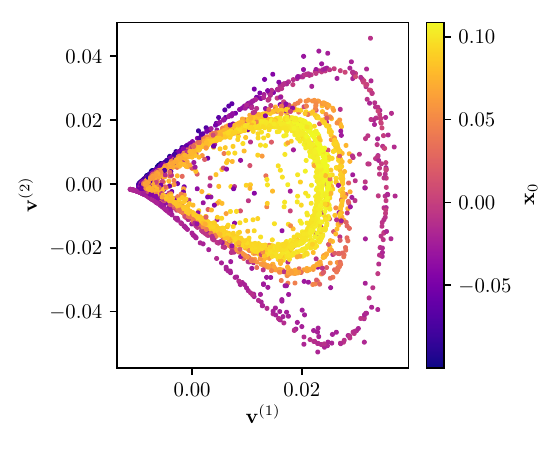}
  \caption{\footnotesize Diffusion maps embedding of the Brusselator system time series windows of length $5$ of the training data. $\boldsymbol{v}^{(1)}$ and $\boldsymbol{v}^{(2)}$ are the two independent diffusion maps modes spanning the data manifold. The color corresponds to one warmed-up internal state variable (here, $\mathbf{x}_0$, one of the 1024 internal reservoir states).}
  \label{fig:brusselator_embedding_space}
\end{figure}

We can now use the diffusion maps-learned starting  map to find the initialization of the reservoir states for any new short input time series window of length $5$. This cold-starting of the reservoir produces more accurate autonomous $H$ steps ahead predictions as we show in Figure~\ref{fig:predictions_new}.
We notice that the autonomous predictions produced by the reservoir which is cold-started with our proposed starting map is much more accurate than the one produced by the traditionally initialized one.

We produce similar experiments with the Lorenz system.
Here, we again sample $600$ trajectories for training, with initial conditions $u_0\sim \mathcal{N}(10,1)$, $v_0\sim \mathcal{N}(1,1)$, and $w_0 \sim \mathcal{N}(0,1)$.

\begin{figure}[h!]
  \centering
  \includegraphics[width=0.7\textwidth]{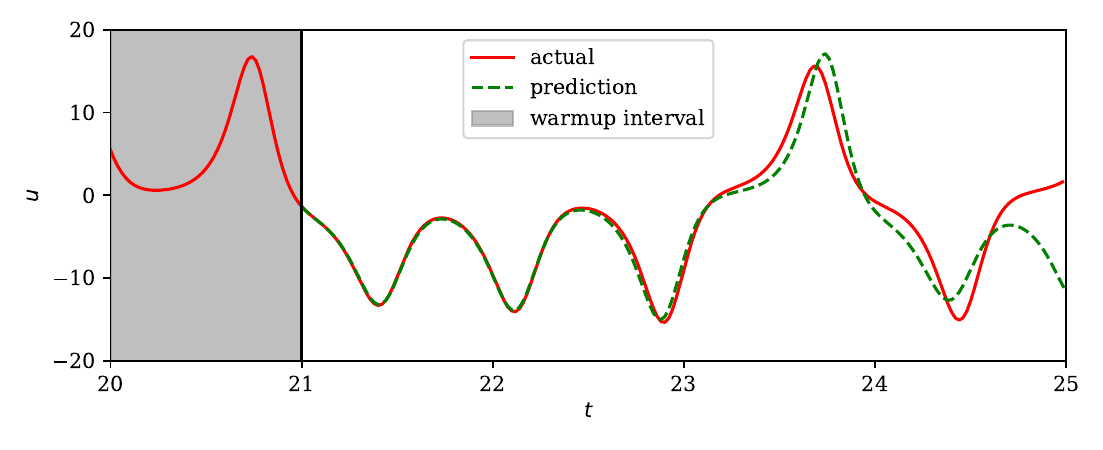}
  \caption{{ \footnotesize Autonomous path-continuing of the partial observations of the Lorenz system produced by the ESN with the readout $\widehat{h}_{ridge}$ and with the initial zero state compared to the true trajectory. The shaded area marks the part of the path which is used as a history to drive the trained ESN; the black line denotes the moment when the subsequent autonomous path-continuation starts.}}
  \label{fig:prediction_lorenz}
\end{figure}

For each trajectory, we sample for $2$ dimensionless time units between $t_{\text{min}}=20$ and $t_{\text{max}}=22$ steps with $\delta t=0.02$, with results in $T=100$ discrete time observations.
Of those trajectories, we discard the first $T_w=20$ ($0.4$ in the intrinsic time of the system) steps washout discretized steps. This results in $600$ pairs of $T_{\text{tr}}=79$ training paths ($1.6$ in the characteristic time of the system).
For testing, we use sample trajectories with $t_{\text{min}}=20$ and $t_{\text{max}}=25$ using $\delta t=0.02$, resulting in trajectories consisting of $250$ time steps (of a duration of $5$ dimensionless time units). One such trajectory is shown in Figure \ref{fig:prediction_lorenz}.
As for the case of the Brusselator system, we create the readout map using ridge regression.
For the reservoir, we use $N=2048$, $\alpha=0.5$, $\lambda = 0.001$, and a spectral radius of $0.80$ for the regression problem. 
Predictions using the thus trained reservoir using a warmup length of $50$ steps are shown in Figure~\ref{fig:prediction_lorenz}. Again, we chose a long warmup length to washout the effect of the reservoir initialization.
 
\begin{figure}[ht]
  \centering
  \includegraphics{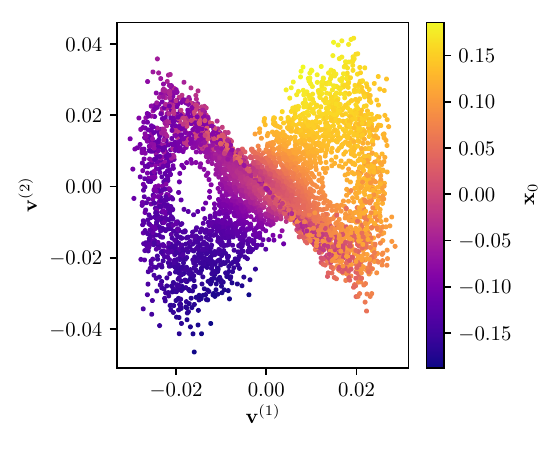}
  \caption{\footnotesize Diffusion maps embedding of the Lorenz system time series windows of length $7$ of the training data. $\boldsymbol{v}^{(1)}$ and $\boldsymbol{v}^{(2)}$ are the first two of the three independent diffusion maps modes spanning the data manifold. The color corresponds to one warmed-up internal state variable (here, $\mathbf{x}_0$, one of the 2048 internal reservoir states). Notice that, though the coloring, this reservoir state is a function of the data manifold.}
  \label{fig:data_manifold_Lorenz}
\end{figure}

Again, we can create a starting map by first approximating the data manifold using diffusion maps and windows of the training time series (here, windows of length $7$).
In this case the data manifold is spanned by three diffusion modes, which is in agreement with the dimension of the original dynamical system. A projection on the first two independent modes is shown in Figure~\ref{fig:data_manifold_Lorenz}.

Finally, we create a mapping from the three-dimensional data manifold to the corresponding internal states of the reservoir. The states are thereby obtained by forcing the reservoir with the training time series, whereas the mapping is again created by fitting geometric harmonics.

\begin{figure}[h!]
  \centering
  \includegraphics[width=0.7\textwidth]{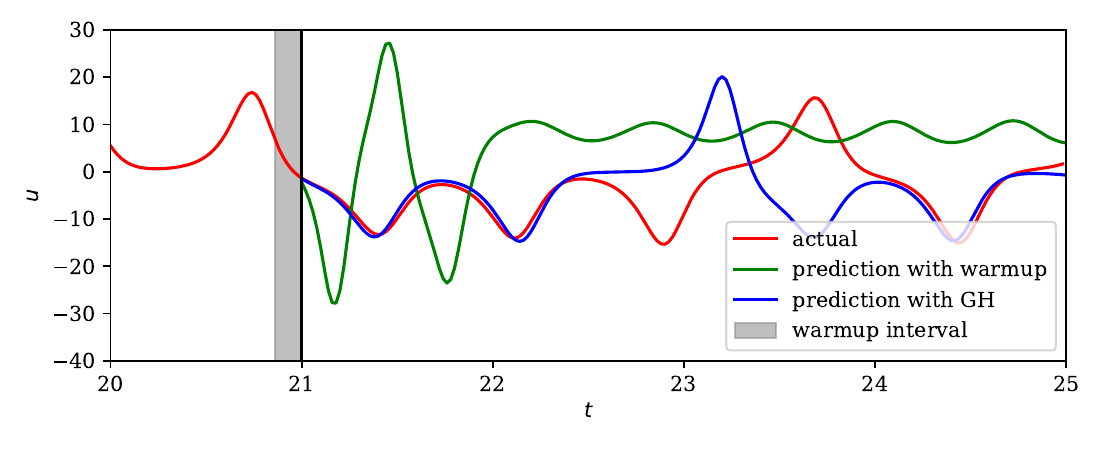}
  \caption{\footnotesize Representative trajectory of the test data for the Lorenz system (red). The warmup period is of length 5 (gray-shaded region). Green - predictions of the ESN with the states initialized as zero vectors and warmup used.
  Blue - predictions of the ESN with the states initialized with geometric harmonics (GH) method.}
  \label{fig:fig_Lorenz_predictions}
\end{figure}

We can now compare the efficacy of our initialization approach versus the traditional warmup approach.
For a short warmup period of just $7$ steps, the prediction results are depicted in Figure~\ref{fig:fig_Lorenz_predictions}. Note that the classical initialization approach leads to a fast divergence of the predicted and true dynamics, since $7$ steps seem to be insufficient to properly warm up the reservoir.
In contrast, using our initialization approach, we obtain forecasts that stay true to the actual dynamics for a long time horizon. Note that due to the approximation errors of the trained resevoir and the starting map, as well as the chaotic nature of the dynamics, predictions will eventually diverge.

\subsection{Robustness in learning the starting map.}
\label{Robustness in learning the cold-start map}
In this section, we empirically study the robustness of our proposed approach with respect to the initialization of the reservoir states using the starting map $\sigma$. More specifically, we explore the sensitivity of the forecasts produced by the reservoir systems with respect to potential imprecisions in the learning of the starting map. As opposed to the previous section where the so-called diffusion maps procedure was used to obtain the starting map $\sigma$ and the linear readout was trained using a ridge regression that admits a closed-form solution, here, instead, we consider other techniques for these steps. In particular, neural network models are employed for approximating $\sigma$ and the reservoir readout $h$. In the following paragraphs, we show that the results obtained with the help of cold-staring initialization of the reservoir systems do not depend much on the particular choice of the learning method. To exemplify this claim, we conduct a series of exercises that consist of the following steps:
\begin{itemize}
    \item Learn the starting map $\sigma$ using a neural network of a given architecture, denote a neural network approximation of $\sigma$ by $\sigma^{\text{NN}}$;
    \item Learn the readout map $h$ using another neural network of a given architecture and denote it as $h^{\text{NN}}$;
    \item Take a set of $10$ arbitrary chosen $2q+1=7$ partial subsequent observations $\boldsymbol{\omega}_k\in \mathbb{R}^{2q+1}$, $k=1,\ldots, 10$ of the Lorenz system to construct the corresponding initial reservoir states $\mathbf{x}_k^{\sigma^{\text{NN}}} = \sigma^{\text{NN}}(\boldsymbol{\omega}_k)$, $k=1,\ldots, 10$;
    \item Construct a set of $1000$ equally distanced values $\sigma_{\eta}^2 \in [0, 0.03]$. For each ${\sigma^j_{\eta}}^2$, $j=1, \ldots, 1000$, a sample of $K=10$ random innovations $\{ \eta^j_k \}_{k\in \{1, \ldots K\}}$, $ \eta_k^j \sim \mathcal{U}\{0, \sqrt{12}\sigma_{\eta}^j\}$, is drawn. 
    \item  Each of the cold-starting states is perturbed  $\tilde{\mathbf{x}}_k^{\sigma^{\text{NN}}} := {\mathbf{x}}_k^{\sigma^{\text{NN}}} + \eta^j_k$, $k=1,\ldots, 10$, $j=1,\ldots, 1000$;
    \item The learnt readout map $h^{\text{NN}}$ is applied to these perturbed initial states and the reservoir system is run autonomously to produce $100$ future steps of the path-continued trajectory;
    \item The mean squared error of the $100$ autonomous predictions is computed per each perturbed state and the corresponding innovation, which results in $10000$ measurements which are subsequently plotted using the scatter plot versus the corresponding values of $\sigma_{\eta}^j$, $j=1, \ldots, 1000$.
\end{itemize}

\begin{figure}[h]
	\centering
	\includegraphics[width=0.5\textwidth]{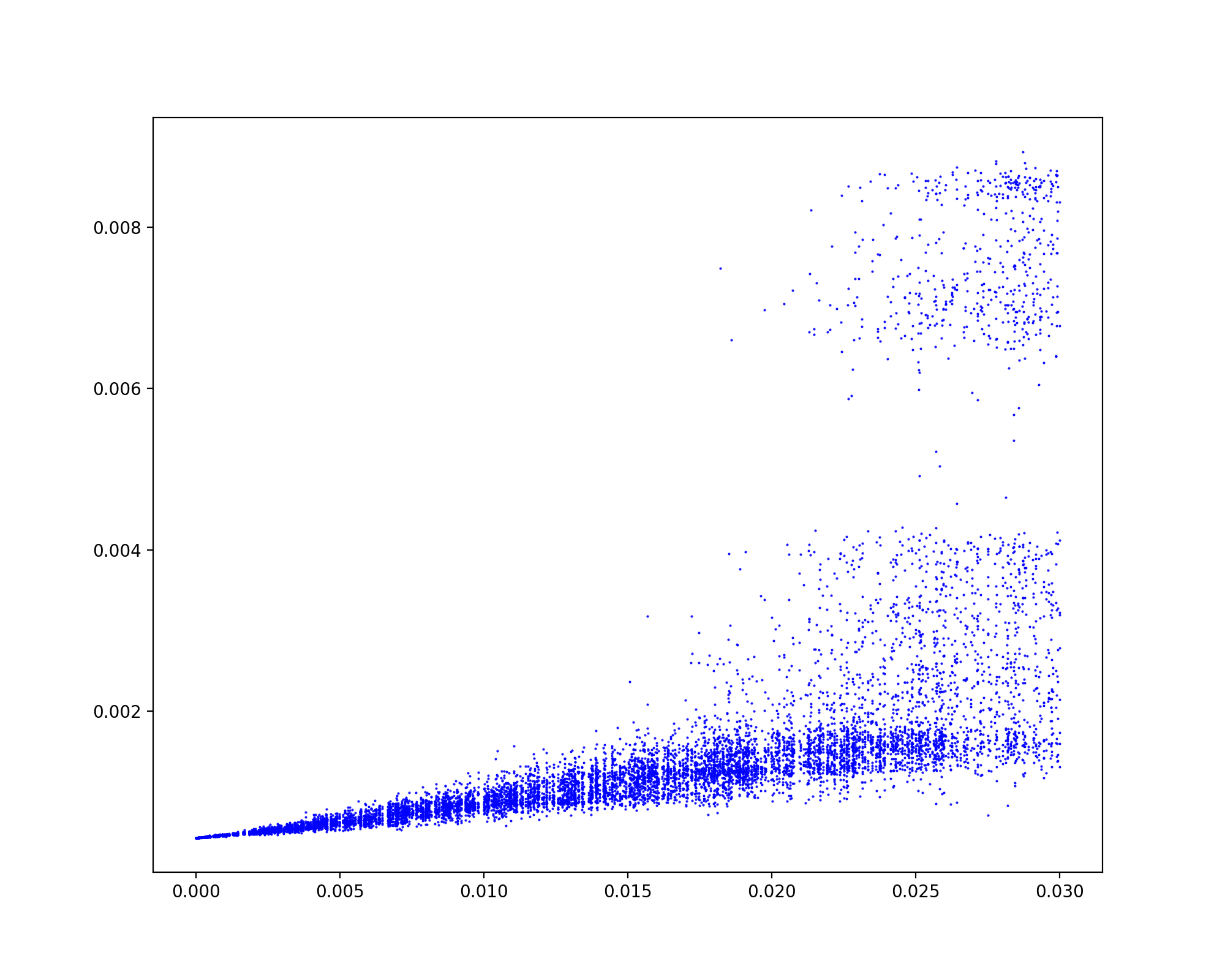}
	\caption{\footnotesize Lorenz system results: mean squared error calculated over 100 time steps versus perturbation of the value suggested by the cold-start map; 10 experiments are conducted for each perturbation error and the perturbation error is varied in the interval 0 to 0.03 with a step-size of $0.03/1000$.}
\label{fig:scatter_lorenz}
\end{figure}

We notice that, in contrast to the previous section, where the readout map was obtained as the solution of ridge regression, the impact of the perturbation on the reservoir outputs is not linear, even in the first step of the out-of-sample prediction. Hence, one expects that the perturbations introduced with respect to the true images of the starting maps get nonlinearly amplified by the neural network readout at the time of autonomous forecasting. Figure~\ref{fig:scatter_lorenz} shows that the dependence of the mean squared forecasting errors as a function of the variance of the perturbing innovations for all the chosen sets of partial subsequent observations $\boldsymbol{\omega}_k\in \mathbb{R}^{2q+1}$, $k=1,\ldots, 10$, is $O(\sigma_\eta^2)$.

The details of our implementation of the learning of the cold-starting and the readout maps are provided in the following paragraphs for the interested reader.

\medskip

\noindent \textbf{Reservoir design.}  We consider the echo state network defined in \eqref{esn equation} with the state dimension $N=900$, and the connectivity (reservoir) matrix $A$ and the input matrix $C$ randomly sampled from $\mathcal{U}\{0,1\}$ distribution. We normalize $A$ such that its spectral radius is $\rho(A) = 0.99$ to satisfy the sufficient condition for the echo state property.  We choose the leak rate to be $\alpha = 0.7$. 

\medskip

\noindent {\bf Learning of the starting map $\sigma$ and the readout $h$ with the neural network approach.} In order to distinguish between the two neural networks used for these two different purposes, we will call them the cold-starting and the readout neural networks, respectively.  To collect the training set for the readout neural network, we use a random initial condition for our discretized Lorenz ODE.  The trajectory of length $T$ of the first coordinate is used as the input of the reservoir system $\{{u_t}\}_{t\in \{1, \ldots, T\}}$ and to collect the corresponding $T$ states of dimension $900$. The washout period of the length $T_{w}$ is discarded and only $T_{\text{tr}} = T-T_w$ pairs of the states and their corresponding one-step ahead values of the trajectory $(\mathbf{x}_t, u_{t+1})$ are used to construct the training set of the length $T_{\text{tr}}$. The neural network is trained with this training set in order to produce the approximation $h^{\text{NN}}$ such that $h^{\text{NN}}(\mathbf{x}_t) =  u_{t+1}$ for all $t$.

In the case of the cold-starting neural network, we proceed as follows. We first observe that 
the co-domain of the starting map $\sigma$ has a  large Euclidean dimension (for example, $N=900$ in our experiment). We hence approximate its image using the $K$ leading principal components ($K$ is arbitrarily chosen; for example we use $K=100$ in our study) obtained by the principal component analysis (PCA). We define a  map $P_{K}:\mathbb{R}^{2q+1} \longrightarrow \mathbb{R}^K$,   so that  
$P_{K}(\omega(\phi^{t-1}(m)), \omega(\phi^{t-2}(m)), \ldots , \omega(\phi^{t-2q-1}(m)))$ 
contains the $K$ leading principal component values of $f(\phi^{t-1}(m))$.
To construct the training set for this neural network we collect the states of the reservoir in the same fashion as above and discard in the same manner the first $T_w$ observations (washout). We compute the projection of the collected data onto the $K$ leading principal components and denote them as $\mathbf{x}^K_t \in \mathbb{R}^K$ for all $t$ in the training set. We use the same notation as above and for every state $\mathbf{x}_t$ denote by $\boldsymbol{\omega}_t:=(\omega(\phi^{t-1}(m)), \omega(\phi^{t-2}(m)), \ldots , \omega(\phi^{t-2q-1}(m)))$ its corresponding history of the inputs-observations. The pairs $(\boldsymbol{\omega}_t, \mathbf{x}^K_t)$ for all $t$ in the training set are used for the cold-starting neural network optimization.  Once the neural network is trained, for any new short history of observations the projection of the corresponding reservoir state onto its $K$ leading principal components is obtained with the neural network. Next, we embed the output of the neural network (the image of the $P_K$ map) into $\mathbb{R}^{N}$ via completing via padding $N-K$ zeroes. Finally, we derive an approximation to the image of $\sigma$  to be the inverse of the PCA transform acting on this zero-padded vector in $\mathbb{R}^{N}$.

\medskip

\noindent {\bf Neural network architectures for the learning of $\sigma$ and $h$.} Throughout our experiments, we use a feedforward neural network used to train 
map $P_{100}$    -- the network is constructed with 4 hidden layers with a
layer dimension equal to 500. 
The activation function on the input and hidden layers is the ReLU
function built into Keras, whereas the output layer has no activation function. Training is
accomplished using the Adam optimizer, minimizing the mean square error as the loss
function. The network is trained using the ReduceLROnPlateau callback function of Keras, which monitors the value of the loss function on the validation set and reduces the learning rate when that loss reaches a plateau. The initial learning rate is set to 0.001, which is halved whenever a plateau of at least 50 epochs is reached. 
While learning $P_{100}$, we
use 500 training epochs and a batch size of 500.
The readout $h^{\text{NN}}$ was also learned with the same 
feedforward network with 5000 state values. While training  $P_{100}$ or $h^{\text{NN}}$, 20\% of the training length was used for validation.


The reader may note that we have not used Lyapunov exponents of the autonomous system resulting from a cold-start to ascertain the robustness of the starting map. This is because the Lyapunov exponents while reflecting on the magnitude of the exponent reflecting the time scale on which system dynamics become unpredictable would depend on the error that would have incurred while learning the readout rather than the error that would have incurred in the cold-start.  
This contrasts the error in the short-term prediction of the learned reservoir with the long (and even infinite) time accuracy of its approximation of the original problem. Given the sensitivity to initial conditions and the lack of guarantees for the smooth dependence of the Lyapunov exponents to small system identification errors, asking for accurate Lyapunov exponent approximation lies beyond the scope of the present work.

\section{Conclusion}
\label{sec:Conclusion}
While observing a solution of an initial value problem with an ordinary differential equation or while iterating a map on an initial condition, one can start observing the solution right away.  The aforementioned amenity was not accessible for forecasting with an echo state network model since even in their autonomous mode, they had to be driven by a 
not-so-short history of the very trajectory that one wanted the model to forecast.  We have overcome this challenge with the notion of a cold start. By employing a small segment of the partial observations (enough to determine a unique state of the underlying dynamical system)
as the initial condition, and using a starting map, we show that it is theoretically possible to initialize the internal state of the reservoir enabling forecasting by iteration from that internal state when the network is run in its autonomous mode. We have also pointed out the natural conditions that entail that the starting map is well-behaved in the sense that it is a Lipschitz function which also justifies the numerically observed robustness of its learning.

From the larger perspective of modeling differential equations, the ``well-trained, well-initialized" reservoir is a numerical approximation of the actual dynamical system. Therefore, the notion of shadowing property would be needed to compare the trajectories of dynamical systems with their numerical approximations~\cite{coven88_pseud_orbit_shadow_famil_tent_maps, grebogi02_shadow_chaot_dynam_system, sauer97_how_long_do_numer_chaot, kennedyil_shadow_higher_dimen}. Some of the authors are currently researching this topic. 

\appendix
\section{Diffusion Maps}
\label{sec:dmaps}

The diffusion maps parametrization technique provides a strategy for the dimensionality reduction of a finite dataset, $X = \{\mathbf{x}_i\}_{i=1}^n$, where each $\mathbf{x}_i \in \mathbb{R}^m$ is a sample from a manifold $M$~\cite{Coifman2006}. 
The first step in the diffusion maps method involves establishing a random walk across the dataset. 
This is facilitated by the creation of an affinity matrix $K \in \mathbb{R}^{n \times n}$, which represents the connections among the points in $X$. The elements of this matrix, $K_{ij}$, are calculated using a kernel, here a Gaussian kernel, according to:

\begin{equation*}
K_{ij} = \exp\left(-\frac{\|\mathbf{x}_i - \mathbf{x}_j\|^2}{2\epsilon}\right),
\end{equation*}
where $\|\cdot\|$ denotes the chosen norm for the data, in this case, the $L_2$ norm. The hyperparameter $\epsilon > 0$ controls the decay rate of the kernel: for smaller values of $\epsilon$, only proximal points are considered connected in $K$, as $K_{ij}$ approaches 0 for distant points.

The diffusion maps algorithm hinges on the normalized graph Laplacian of the data converging to the Laplace-Beltrami operator on the manifold $M$ as the number of points $n \rightarrow \infty$ and $\epsilon \rightarrow 0$. However, a specific normalization is required for data obtained from non-uniformly sampled points to accurately recover the Laplace-Beltrami operator. This involves defining a diagonal matrix $D \in \mathbb{R}^{n \times n}$, with $D_{ii}=\sum_{j=1}^{n} K_{ij}$, and then calculating the normalized affinity matrix, given by
\begin{equation*}
\tilde{K} = D^{-\kappa}KD^{-\kappa},
\end{equation*}
where $\kappa$ modulates the density effect. For $\kappa = 0$, the density's influence is maximal, suitable only for uniformly sampled data, whereas $\kappa = 1$ removes the density effect, enabling the recovery of the Laplace-Beltrami operator~\cite{coifman08_graph_laplac_tomog_from_unknow_random_projec}. Another normalization step yields $S$, a Markovian matrix, by dividing each entry of $\tilde{K}$ by the sum of its rows. The eigendecomposition of $S$ reveals a complete set of real eigenvectors $\boldsymbol{v}^{(i)}$ and eigenvalues $\lambda_i$, facilitating a nonlinear parametrization of the dataset $X$ in terms of these eigenvectors. Selecting the leading eigenvectors that are independent/non-harmonic generates a set of latent variables $\Phi = \{\boldsymbol{v}^{(1)}, \ldots , \boldsymbol{v}^{(d)}\}$ that encapsulate the intrinsic geometry of the manifold from which the dataset was sampled. If the number of these selected eigenvectors $d$ is less than the original variable dimensions $m$, the process effectively reduces dimensionality by presenting a more simplified representation of the dataset. For a dataset $X$ comprising short time series windows $u_t$, diffusion maps enable the extraction of reduced latent variables in a data-driven manner, with $\kappa = 0$ and $\epsilon$ chosen as the median of all pairwise distances, ensuring that the choice of $\alpha$ does not qualitatively alter the diffusion map results.

\section{Geometric Harmonics}
Geometric harmonics is utilized to extend a function $\mathcal{F}$, potentially vector-valued, sampled at certain points \(X=\{\mathbf{x}_i\}\) on a manifold \(M\), to a new point \(\mathbf{x}_{\text{new}} \notin X\)~\cite{Coifman2006}. 
In this context, a modified approach of geometric harmonics is employed to interpolate \(\mathcal{F}\) using the reduced coordinates \(\Phi\) identified through diffusion maps. Specifically, after the dimensionality reduction phase yields non-harmonic eigenvectors, the goal is to express \(\mathcal{F}\) in terms of these reduced coordinates \(B = (\boldsymbol{v}^{(1)}|\boldsymbol{v}^{(2)}|\ldots|\boldsymbol{v}^{(d)})\in \mathbb{R}^{n\times d}\) with $\boldsymbol{v}^{(j)} \in \mathbb{R}^n$. Despite the exclusion of harmonic eigenvectors, a subsequent application of diffusion maps to the coordinates \(\Phi\) facilitates the creation of a functional basis connecting \(\Phi\) to any function \(\mathcal{F}\) defined on the original space.

Similar to the initial diffusion maps process, the first step involves calculating an affinity matrix \(C_{i,j} = C(\boldsymbol{b}_{i}, \boldsymbol{b}_{j}) = \exp\left(-\frac{\|\boldsymbol{b}_{i} - \boldsymbol{b}_{j}\|_2^2}{2\epsilon'}\right)\), where $\boldsymbol{b}_i \in \mathbb{R}^d$ denotes the $i$-th row of the matrix $B$. Being symmetric and positive semidefinite, \(C\) possesses orthonormal vectors \(\boldsymbol{\psi}^{(1)}, \boldsymbol{\psi}^{(2)}, \ldots, \boldsymbol{\psi}^{(n)}\) and non-negative eigenvalues \(\sigma_1 \geq \sigma_2 \geq ... \geq \sigma_n \geq 0\). These eigenvectors serve as a projection basis for extending a function \(\mathcal{F}\). Selecting a threshold \(\delta > 0\), the set of significant eigenvalues \(S_{\delta} = \{\alpha : \sigma_{\alpha} > \delta\sigma_1\}\) is determined, where \(\delta\) is chosen such that \(d < {\rm{Card}}(S_{\delta}) < n\). 

Projecting the image of \(\mathcal{F}\) onto this truncated eigenvector set yields an approximation \(\mathcal{F} \approx P_{\delta}\mathcal{F} \equiv \tilde{\mathcal{F}} = \sum_{\alpha \in S_{\delta}} \boldsymbol{\psi}^{(\alpha)} (\tilde{\mathcal{F}}^T \boldsymbol{\psi}^{(\alpha)})^\top\).

To extend \(\tilde{\mathcal{F}}\) to a new coordinate \(\boldsymbol{b}_{\text{new}}\), which is not one of the rows of $B$, the extension is given by \(\tilde{\mathcal{F}}_{\text{new}}(\boldsymbol{b}_{\text{new}}) = \sum_{\alpha \in S_{\delta}} \boldsymbol{\psi}^{(\alpha)}_{\text{new}}  (\tilde{\mathcal{F}}^\top  \boldsymbol{\psi}^{(\alpha)})^\top\), with \(\boldsymbol{\psi}^{(\alpha)}_{\text{new}} = \sigma^{-1}_{\alpha} \sum_{i=1}^{n} C(\boldsymbol{b}_{\text{new}}, \boldsymbol{b}_{i}) \cdot \psi^{(\alpha)}_i\) and where \(\psi^{(\alpha)}_i\) is the \(i\)-th component of the eigenvector \(\boldsymbol{\psi}^{(\alpha)}\). This approach, employing a truncated set \(S_{\delta}\), addresses numerical instabilities that occur when \(\sigma_{\alpha} \rightarrow 0\).

By applying geometric harmonics in this manner, it is possible to predict the values of \(\mathcal{F} = \mathbf{x}_t\) at unseen points \(\boldsymbol{b}_{\text{new}}  \in \mathbb{R}^d\), for $d=2$ for the Brusselator system, $d=3$ for the Lorenz system, is derived via Nyström extension~\cite{nystroem30_ueber_prakt_aufloes_von_integ} on time series windows of \(u_t\).

\medskip
\addcontentsline{toc}{section}{Supplementary information}
\noindent {\bf Supplementary information.}
All code necessary to reproduce the numerical results presented in the paper are publicly available at \url{https://github.com/Learning-of-Dynamic-Processes/coldstart}.

\medskip

\addcontentsline{toc}{section}{Acknowledgments}
\noindent {\bf Acknowledgments.} LG and GM thank the hospitality of Nanyang Technological University, where part of this work was completed. GM acknowledges partial funding through an incentive grant, UID 150668 from the NRF, South Africa. JPO acknowledges partial financial support from the School of Physical and Mathematical Sciences of the Nanyang Technological University. BH acknowledges financial support from   the Air Force Office of Scientific Research under MURI award number FA9550-20-1-0358 (Machine Learning and Physics-Based Modeling and Simulation) and the Department of Energy under the MMICCs SEA-CROGS award. The work of YK and FK was partially supported by DARPA and the US Air Force Office of Scientific Research. 


\addcontentsline{toc}{section}{Bibliography}
\bibliographystyle{wmaainf}
\bibliography{GOLibrary, lit}

\end{document}